\documentclass[11pt]{amsart}
\usepackage{amsxtra}
\usepackage{amssymb}
\usepackage{graphicx}

\newenvironment{dedication}
{
	\thispagestyle{empty}
	\vspace*{\stretch{1}}
	\itshape             
	\raggedleft          
}
{\par 
	\vspace{\stretch{3}} 
}

\theoremstyle{plain}

\numberwithin{equation}{section}

\newmuskip\pFqmuskip

\newcommand*\pFq[6][8]{%
	\begingroup 
	\pFqmuskip=#1mu\relax
	\mathcode`\,=\string"8000
	\begingroup\lccode`\~=`\,
	\lowercase{\endgroup\let~}\pFqcomma
	{}_{#2}F_{#3}{\left[\genfrac..{0pt}{}{#4}{#5};#6\right]}%
	\endgroup
}
\newcommand{\pFqcomma}{\mskip\pFqmuskip}

\newtheorem{theorem}{Theorem}[section]
\newtheorem{lemma}[theorem]{Lemma}
\newtheorem{proposition}[theorem]{Proposition}
\newtheorem{corollary}[theorem]{Corollary}
\newtheorem{definition}[theorem]{Definition}
\newtheorem{example}[theorem]{Example}

\newtheorem{remark}[theorem]{Remark}
\newtheorem{conjecture}[theorem]{Conjecture}

\newtheorem{question}[theorem]{Question}
\newcommand \bth[1] { \begin{theorem}\label{t#1} }
	\newcommand \ble[1] { \begin{lemma}\label{l#1} }
			\newcommand {\ethe} { \end{theorem} }
		\newcommand {\ele} { \end{lemma} }
		\newcommand \bpr[1] { \begin{proposition}\label{p#1} }
			\newcommand \bco[1] { \begin{corollary}\label{c#1} }
		\newcommand \bde[1] { \begin{definition}\label{d#1}\rm }
					\newcommand \bex[1] { \begin{example}\label{e#1}\rm }
						\newcommand \bre[1] { \begin{remark}\label{r#1}\rm }
							\newcommand \bcon[1] { \begin{conjecture}\label{con#1}\rm }
								\newcommand \bque[1] { \begin{question}\label{que#1}\rm }
			\newcommand{\beq}{\begin{equation}}
			\newcommand{\eeq}{\end{equation}}
			\newcommand{\beqa}{\begin{eqnarray}}
			\newcommand{\eeqa}{\end{eqnarray}}
			\newcommand{\beaa}{\begin{eqnarray*}}
				\newcommand{\eaa}{\end{eqnarray*}}

		\newcommand {\epr} { \end{proposition} }
						\newcommand {\eco} { \end{corollary} }
					\newcommand {\ede} { \end{definition} }
				\newcommand {\eex} { \end{example} }
			\newcommand {\ere} { \end{remark} }
		\newcommand {\econ} { \end{conjecture} }
	\newcommand {\eque} { \end{question} }

\newcommand \thref[1]{Theorem \ref{t#1}}
\newcommand \leref[1]{Lemma \ref{l#1}}
\newcommand \prref[1]{Proposition \ref{p#1}}

\newcommand \exref[1]{Example \ref{e#1}}

\def \A {{\mathcal A}}
\def \B {{\mathcal B}}

\def \K {{\mathcal K}}
\def \M {{\mathcal M}}
\def \R {{\mathcal R}}
\def \F {{\mathcal F}}

\def \Cset {{\mathbb C}}

\def \Nset {{\mathbb N}}

\def \span { {\mathrm{span}} }

\def \ad { {\mathrm{ad}} }

\newcommand{\al}{\alpha}
\newcommand{\be}{\beta}

\def \Cset {{\mathbb C}}

\def \Nset {{\mathbb N}}

\def \B  {{\mathcal{B}}}               

\def \al {\alpha}
\def \be {\beta}

\def \la {\lambda}
\def \La {\Lambda}

\def \ga {\gamma}

\def \Ga {\Gamma}

\def \n  {\mathfrak{n}}

\newmuskip\pFqmuskip

\usepackage{hyperref}

\begin{document}

\title[ VOP with Bochner's property]
 {Vector orthogonal  polynomials with Bochner's property
  }

 	\author[E.~Horozov]{Emil Horozov}

 	\address{
 			Department of Mathematics and Informatics \\
 			Sofia University \\
 			5 J. Bourchier Blvd. \\
 			Sofia 1126 \\
 			Bulgaria, and\\	
 			Institute of Mathematics and Informatics, \\ 
 			Bulg. Acad. of Sci., Acad. G. Bonchev Str., Block 8, 1113 Sofia,
 			Bulgaria	}
 	\email{horozov@fmi.uni-sofia.bg}

\date{\today}

\date{}

\begin{abstract}

	Classical orthogonal polynomial systems of Jacobi, Hermite, Laguerre and Bessel have the property that the polynomials of each system are eigenfunctions of a second order ordinary differential operator. According to a classical theorem by Bochner they are the only systems  with this property.  Similarly, the polynomials of Charlier, Meixner, Kravchuk and Hahn are both orthogonal and are eigenfunctions of a suitable difference operator of second order. We  recall that according to the famous theorem of Favard-Shohat, the condition of orthogonality is equivalent to the 3-term recurrence relation.  	Vector orthogonal polynomials  (VOP) satisfy  finite-term recurrence relation with more terms,  according to a theorem by J. Van Iseghem  and this characterizes them.
	
 Motivated by Bochner's theorem  we are looking for VOP that are also eigenfunctions of a differential (difference) operator. 
	We call these simultaneous conditions Bochner's property.   
 
	The  goal of this paper is to introduce methods for construction of  VOP which  have Bochner's property.  The methods are purely algebraic and are based on automorphisms of non-commutative algebras. They also use ideas from the so called bispectral problem.   Applications of the abstract methods include  broad generalizations of the classical  orthogonal polynomials, both continuous and discrete. Other results connect different families of VOP, including the classical ones, by  linear transforms of purely algebraic origin, despite of the fact that, when interpreted analytically,  they  are integral transformations.   
\end{abstract}

\maketitle

\begin{dedication}
	To Boris Shapiro on occasion of  his 60-th birthday
\end{dedication}

\medskip

\section{Introduction}  \label{intro} 

The  present paper generalizes the results of \cite{Ho1, Ho2}. There we sketched a program aiming to construct sets of  Vector Orthogonal Polynomials (VOP), with the property that the polynomials of each system  are eigenfunctions of some differential or difference operator. The problem stems from   the famous Bochner theorem. Let us recall it.

Bochner   \cite{Bo} has proved that all systems of orthogonal polynomials $P_n(x)$, $n=0, 1, \ldots$, $\deg P_n = n$      that are also  eigenfunctions of a second order differential operator

\begin{equation}
L(x, \partial_x) = A(x)\partial_x^2 + B(x)\partial_x + C(x)  \label{BP1}
\end{equation}
with eigenvalues $\lambda(n)$, are exhausted by  the classical orthogonal polynomials of Hermite, Laguerre and Jacobi and the lesser known Bessel ones  (found in \cite{KF}).
 
The orthogonality condition, due to the classical theorem by Favard and Shohat is equivalent to  the well known 3-term recurrence relation

\begin{equation*}
xP_n = P_{n+1} + \beta(n)P_n + \gamma(n) P_{n-1},   
\end{equation*}
where $\beta(n), \gamma(n)$ are constants, depending on $n$.\footnote{Here  and further in the paper  the polynomials $P_n(x), \; \; n=0, 1, \ldots,$    are
	monic, i. e. they are normalized by the condition that the coefficient at $x^n$ is $1$.}  Also throughout this paper we use orthogonality in slightly broader sense - with respect to some functional, which  does not need to be defined by a measure, nor to be positive-definite.

 For example  the   Hermite polynomials $\{H_n(x),\; n=0,1, \ldots \}$ are  eigenfunctions of the Hermite differential operator:

\[
L = -\partial_x^2 + x\partial_x
\]
with eigenvalues $\la(n)= n$. At the same time they satisfy 3-term recurrence relation:

\[
xH_n(x) = H_{n+1} + nH_{n-1}.
\]

\bre{Her}
The version of Hermite polynomials here is called   "probabilists' Hermite polynomials". However in literature the so called "physicists' Hermite polynomials" are more widely used,  see e.g. \cite{NIST}. Naturally they also can be obtained by our methods and in fact both versions are rescaling of one another. Here we use the probabilists' version as more suitable for the methods.
\ere

In 1940 H. Krall \cite{Kr} extended Bochners's result to include the case of differential operators of order 4. In recent times there is much activity in generalizations and versions of the classical results of Bochner and Krall. In particular some necessary conditions  for  a polynomial system to consist of eigenfunctions of a differential operator can  be found in \cite{EKLW} and in the references therein.

An important  role in some of these generalizations plays the ideology of the bispectral problem which was initiated in \cite{DG}. Translating the Bochner's and Krall's results into this language  presents a good basis to continue investigations, see \cite{GH3, GH4}. It goes as follows. Let us introduce the bivariate function  $\psi(x, n) = P_n(x)$ and the shift operators $T$ and $T^{-1}$, acting on the discrete variable $n$.  If we write the right-hand side of the 3-term recurrence relation as a difference operator $J(n)$ acting on  $n$,  then the  3-term recurrence relation can be written as

\[
J(n) \psi(x, n)  = x\psi(x, n). 
\]
where in terms of the shift operators $T^{\pm}$, which map $f(n)$ into  $f(n\pm1)$   the operator $J(n)$ has the form:

\beq
J(n) = T + \beta(n)T^0 + \gamma(n)T^{-1}   \label{Diff}
\eeq
This means that $\psi(x, n)$ is an eigenfunction of the discrete operator $J$ with  eigenvalue $x$ and of a differential operator  $L(x, \partial_x)$ with eigenvalue 
$\lambda_n$.

 Hence we can formulate the Bochner-Krall problem as the following bispectral problem:
 
  \begin{quote}
  	
  	(BK)
   \textit{
 	Find all systems of orthogonal  polynomials $P_n(x)$  which are eigenfunctions of a differential operator $L$ \eqref{BP1}   and a difference operator $J$  \eqref{Diff}.
 }
  \end{quote}

   The same  problem for discrete orthogonal polynomials was solved by O. Lancaster \cite{Lan} and P. Lesky \cite{Les}, although earlier E. Hildebrandt \cite{Hil} has found all needed components of the proof. For more information see the excellent review article by W. Al-Salam \cite{AlS}.  In the rest of the introduction we will discuss  mostly the continuous problem as the discrete
    one is very similar.

 The implementation of ideas from the bispectral  problem in these studies is  usually done through the well developed  machinery of Darboux transformations in different versions.  See \cite{MS} for the general theory and \cite{GH, GHH, GY} for applications to orthogonal polynomials. Other methods for solving the bispectral problem have been developed in  \cite{Du, DdlI1, DdlI}.  One starts with some known polynomial sequences and performs Darboux transformations in order to obtain new  systems  of polynomials that are eigenfunctions of a differential operator and at the same time are eigenfunctions of a difference operator $J(n)$ as in \eqref{Diff}. See, e.g \cite{GH3,  GH4, GHH, GY, Il1}.

Here we  also exploit some ideas from the theory of bispectral operators.    Before explaining them and the main results let us introduce one more concept which is central for the present paper. This is the notion of vector orthogonal polynomials  (VOP), introduced by J. van Iseghem \cite{VIs}.   Let $\{P_n(x),\;  n= 0, 1, \ldots\}$ be a family of monic polynomials such that $\deg P_n = n$. The polynomials are VOP iff   there exist $d$  functionals $\mathcal{L}_j, j =0,\ldots,d-1$    on the space of all polynomials  $\Cset[x]$   such that

\[
\begin{cases}
\mathcal{L}_j (P_nP_m)    = 0, m > nd+ j, n \geq 0,  \\
\mathcal{L}_j(P_nP_{nd+ j})    \neq 0, n   \geq 0,    
\end{cases}
\]
for each $j \in N_{d+1} := \{0, \ldots, d-1 \}$. When $d = 1$ this is the ordinary notion of orthogonal polynomials.  The orthogonality connected with    $d$ functionals rather than with only one, gives the other name of VOP - \textit{d-orthogonal polynomials}, which is quite popular.   By  a theorem of   J. Van Iseghem \cite{VIs} (see also  P. Maroni \cite{Ma}), which generalizes the result of Favard and Shohat,  this property is equivalent to the following one: the polynomials $P_n(x)$   satisfy a $d+2$-term recurrence relation 

\beq  \label{d-ort}
xP_n(x) = P_{n+1}(x) + \sum_{j=0}^{d}\gamma_j (n)P_{n-j}(x)
\eeq
with constants $\gamma_j(n)$, independent of $x$, $\ga_d(n) \neq 0, n \geq d$.

The VOP and the broader family of multiple orthogonal polynomials are subject to an intensive research in last sveral decades due to their applications to Hermite-Pad\'e approximation, random matrices, number theory, etc. For a detailed account on the subject we refer to many excellent papers, including review ones \cite{ApKu, Apt, ACVA, BDK, DBr}.
In what follows, speaking about VOP, we always will understand polynomials systems with the above $d+2$-term recurrence relation \eqref{d-ort}.  

We can state the generalized Bochner problem (GBP) as follows.

\begin{quotation}

\textit{Find systems of polynomials $P_n(x), \; n=0, 1, \ldots$,    $\deg P_n = n$, that are eigenfunctions of a 
	differential/difference  operator $L$ of order $m$ with eigenvalues $\lambda(n)$ depending on the discrete variable $n$ (the index):}

\[
LP_n(x) = \lambda(n)P_n(x)  
\] 
\textit{and which at the same time are eigenfunctions of a difference operator, i.e. they satisfy a finite-term (of fixed length $d+2$), recurrence relation of the form} 
\[
xP_n(x) = P_{n+1}(x) + \sum_{j=0}^{d}\gamma_j(n)P_{n-j}(x).
\]

\textit{For each $m$ and each $d$  classify all systems of $d$-orthogonal polynomials.}

\end{quotation}

\vspace{0.3cm}

We will call the polynomial systems described above "Polynomial systems with generalized  Bochner' property". 
Obviously, the problem can be stated for difference, $q$-difference and other types of operators. In this paper we will deal   with both differential and difference operators.

The first part of the GBP has been addressed in many works, e.g. \cite{Dou, BCD2, VAC, vAssc} and references therein. The second part of the problem looks much more difficult. It is solved only for $d=1$, i.e. for orthogonal polynomials and differential operators of order 2 (Bochner) and order 4 (Krall). However some other  cases seem to be in reach of the methods  at our disposal.

It is    evident that we still need to accumulate more examples of polynomial systems with the Bochner's property. If possible, their construction should be done in a systematic way, that eventually will prompt the limitations of the methods. But any individual example would also help to understand better the problem. 



The    goals of the present paper are several. The first one is to   obtain families of VOP with Bochner's property. These families include the classical OP of  Hermite and Laguerre, as well as the discrete systems of Charlier, Kravchuk and Meixner  and their quite straightforward extensions. The next goal is to show the connections between continuous and discrete VOP.  We point out that these  two classes  of VOP             are   mapped into each other in an almost trivial way.  Further we give a number of examples, both old and new, that deserve attention even by themselves. Finally but probably most important is   our general framework which allows to view all such systems as simple realizations of appropriate pairs of associative noncommutative algebras and mappings between them. 
It is worth mentioning  that  all the polynomial systems in this paper originate from the simplest one   $\{x^n, \;  n=0, 1, \ldots \}$. 


We will explain our construction  on the  example  of Hermite-like polynomials.  We start with the   algebra $W_1$ of differential operators in one variable $x$ with polynomial coefficients, which is known under the name the first  Weyl algebra. We define the operator $H= x\partial_x$ and the system of polynomials $\{\psi(x, n) = x^n, \; n=0,  1, \ldots \}$.   Notice that

\[
H\psi(x,n) = n\psi(x,n),
\]
i.e. the functions $\psi(x,n)$ are eigenfunctions of $H$ for each fixed $n$.   The function  $\psi(x,n)$ is also an eigenfunction of the shift operator $T$, but considered as a function in $n$ :

\[
T \psi(x,n) = x    \psi(x,n),
\]
with eigenvalue $x$.
This is a very simple example of the GBP. From it we will obtain  our VOP as follows.
Together with the Weyl algebra we consider the algebra $\R$ generated by the operators $T, T^{-1}$ and the operator of multiplication by $n$. Next we 
  define an anti-homomorphism    $b\colon W_1 \to \R$ on the generators of $W_1$ by 

\[
b(x) = T, \; b(\partial_x) = nT^{-1}.
\]
 From this definitions  it follows that
 
 \[
 b(H) = n.
 \]
 For later use put   $\B_1 =W_1$. Denote by $\B_2$ the image of $\B_1$ under the map $b$. 
 Using  $\psi(x,n) = x^n$    we can express the map $b$ as 

\[
A \psi(x,n) = b(A)\psi(x,n),   \; A \in W_1.
\]
 It is easy to see that $b\colon \B_1 \to \B_2$ is an anti-isomorphism.
 
Now we can perform an automorphism of the Weyl algebra. We notice that all of the automorphisms of $W_1$ were described by Dixmier in \cite{Dix}. In terms of the above operators $x, \partial_x$ all the automorphisms are generated by 

\[
e^{\ad_{G}}, \; \text{with} \; G= \la \partial_x^q, \; \text{or} \; G= \mu x^p \; p, q \in \Nset,\; \la, \mu \in \Cset,
\]
where $\ad_A(B) = [A, B]$, $A, B \in W_1$.
One can easily show that the above operators are well defined as    both $x$ and $\partial_x$ act locally nilpotently, i.e. when applied to any $A \in W_1$ the exponentials defining the above series are finite sums.

 We will take one of the simplest automorphisms $\sigma_G = e^{\ad_G}$, defined via $G= -\partial_x^2/2$. Short computations show that

\[
\sigma_G(H) = -\partial_x^2 + x\partial_x. 
\]
This is exactly the Hermite operator. Let us map the polynomial system $x^n$ into  

\[
P_n(x) = e^G (x^n) = \sum\limits_{j=0}^{\infty}\frac{G^j (x^n)}{j!}.
\]
This is again a polynomial system as the series is  finite and $\deg P_n = n$. If we present $e^{\ad_G}(H)$ as

\[
e^{\ad_G}(H) = e^G H e^{-G},
\]
we easily see that 

\[
\sigma_G(H)P_n(x) =  e^G H e^{-G}P_n(x) = e^G H x^n = n P_n(x),
\]
i.e. the  new polynomials are eigenfunctions of the operator $\sigma_G(H)$.
Notice that the anti-automorphism $b$ is transformed into $b' = b \circ \sigma^{-1}$. It   is easy to show that $b'(x)= T + nT^{-1}$, i.e. it gives the 3-term recurrence of Hermite polynomials:

\[
xP_n(x) = P_{n+1}(x) + nP_{n-1}(x).
\]

\bigskip

We can adapt this construction by: 1)  taking an appropriate pair of algebras $\B_1$, $\B_2$, instead of $W_1$ and $b(W_1)$ and a map $b$ between them and 2) performing automorphisms $\sigma\colon \B_1 \to \B_1$, followed by constructing new anti-automorphism $b' = b\circ \sigma^{-1}$. It turns out that the construction can be performed in a very general situation as follows.


Consider the  Weyl algebra $W_1$ over a field $k$. For simplicity we consider that the field $k$ has characteristic $0$. We recall that $W_1$   is  a quotient of the free algebra on two generators, Z and Y, by the two-sided ideal generated by the element    $ZY - YZ -  1$. We are going to construct a pair of
associative unital algebras  $\B_1, \B_2$. Let $H= YZ$. For  a polynomial $R(H)$   in $H$ put $G  =  R(H)Z$  and define the algebra $\B_1$   

\[
\B_1 = \span{(Y, H, G)}. 
\]
 One can check that $e^{\ad_G}$ acts locally nilpotently on $\B_1$.

 Let us define another algebra $\R$ over $k$, defined by  generators $T, T^{-1}, \hat{n}$ and relations 
 
 \[
 T\cdot T^{-1} =  T^{-1}\cdot T =    1, \; [T, \hat{n}] = T,\; [T^{-1}, \hat{n}] = - \hat{n}.
 \] 
We define an anti-homomorphism $b\colon \B_1 \to \R$ by 

\[
\begin{cases}
b(Y)= T\\
b(H) = \hat{n}\\
b(G) = \hat{n}T^{-1}R(\hat{n}).
\end{cases}
\]
The algebra $\B_2$ will be defined as the image of  $\B_1$ via $b$.  We note that instead of $G$ we can use any polynomial $q(G)$. Of particular interest are the polynomials $q(G) = \rho G^l$ due to the fact that they have more properties,  see \cite{Ho3}.
Similarly to the construction of the Hermite polynomials, corresponding to  $G= \partial_x$ and $q(G) = -\partial_x^2/2$ we define automorphisms
$\sigma_q\colon \B_1 \to \B_1$ and  new anti-isomorphisms $b' = b\circ \sigma_q^{-1}$, which encode all the properties of the polynomial system, that  we even haven't yet defined.  It is easy to define   a polynomial system in each concrete case but also abstractly, which we don't need.   

We notice that this construction while mimicking the case  when the Weyl algebra is realized via differential operators, works perfectly in the case of difference operators, too.  We need this in order to extend some of the classical discrete orthogonal polynomials. The point is that the Weyl algebra can be realized also by   difference operators.   

Our next step is to realize the abstract construction in cases of VOP. The  continuous polynomial systems can  be obtained by using the standard representation of  the Weyl algebra as the algebra of differential operators with polynomial coefficients in one variable $x$. Namely we put 
 $Y$ to be the operator of multiplication by $x$ and $Z$ to be the differentiation $\partial_x$. In this way we obtain extensions of the Hermite and the Laguerre polynomials, similar to those obtained in \cite{Ho1, Dou, BCD2, BCO}. 

The  discrete polynomials are treated similarly but instead of the polynomial system $x^n$ we take the falling factorials $(x)_n = x(x-1)\ldots(x-n+1)$. Then we can realize the Weyl algebra by the
difference operators $D, D^{-1}$,  acting on functions of the variable $x$ as shift operators  $D^{\pm}f(x) = f(x\pm1)$, together with the operator $\hat{x}$, acting as multiplication by $x$. Put $\Delta = D-1$, $\nabla = D^{-1} - 1$, $H = -\hat{x}\nabla$. It is easy to check that  $[\Delta, \hat{x}D^{-1}] = 1$.  In this way we obtain VOP, which  extend the results of \cite{Ho2, BCO, GVZ} on    Charlier and Meixner type polynomials. 

Naturally having   similar constructions we can ask if there is a connection between the parallel continuous and discrete families. We obtain the mapping of these pairs in a very simple way - by the Mellin transform. We have to notice that even in the case of the mapping between Laguerre and Meixner polynomials our construction is new.    A connection of this type  has been  found in  \cite{Koo, AA}. Our construction gives another    formula, which is both very simple and very easy to prove (modulo the general construction).

  The success of this program depends on the possibility to find pairs of algebras   $\B_1, \B_2$ and automorphisms of $\B_1$  such that  the functions defined by the latter series are polynomials. We notice that there are automorphisms that do not share the last property, e.g. the one defined by $G=\hat{x}$. Similar difficulty is not met in cases of differential-differential bispectral problems, cf. \cite{DG, BHY}.
  
  We  note that there are   numerous nice constructions of systems of polynomials, both the classical ones or multiple orthogonal polynomials and even more general special functions using operational methods. One feature of our approach is that we construct not only the polynomial systems but simultaneously many of the accompanying objects like finite-term relations, differential operators, raising and lowering operators, etc. In the present paper we restrict ourselves only with properties that come  automatically with the construction, described above. Other properties of the VOP with Bochner property will be considered in \cite{Ho3}, which is a continuation of this  paper. In particular we will find hypergeometric representations, which will point to  the similarity   of the systems of VOP,  corresponding to    $q(G) = \rho G^l$  with the Gould-Hopper polynomials \cite{GH}. 

During the years after their discovery Gould-Hopper polynomials have received a lot of attention from scientists in various fields. Needless to say that the majority of these are  in special functions, see e.g.    etc., \cite{DSZ,  DTC, SM, DLMTC} and the references therein. However  it is worth noticing that specialists in other fields have also used them. I will mention in particular      combinatorics \cite{Rio, VL}, integrable systems \cite{Cha}. Also connections to quantum physics and probability are found.

The paper is divided into three parts. We start with a brief review in Sect. \ref{BT} of well known material from the bispectral problem, taken mainly from \cite{BHY}.  Sect.  \ref{AC} is the core of the paper. Here we perform our abstract construction, which  consists of description of certain pairs of associative algebras and different maps between them in the spirit of the above explanations. These maps in fact encode all the information about eventual  VOP, the polynomials themselves being left to  Part \ref{II}. In Part    \ref{II}  we  apply the abstract construction,    first in Sect. \ref{con} to differential GBP. This means that we list pairs of algebras of $\B_1, \B_2$ such that $\B_1$ is realized by  differential operators in one variable. Although some  families are well known (from \cite{BCD2} and \cite{Ho1}), most of the VOP are new. In the next Sect. \ref{dis}. we do the same for a realization of $\B_1, \B_2$ in terms of difference operators. Here the new families are almost all. Only the extension of the Charlier  (\cite{BCZ, VZh} and Meixner polynomials \cite{Ho2} are known.  

In the next Sect. \ref{DdsC} we point out explicit maps between the constructed continuous and discrete families via the Mellin transform. 

The last Part \ref{III} is devoted to examples. Many of them  are well known and serve only to illustrate our methods in Sect.  \ref{ee}. There is a second class of known examples, which then are extended to new ones in  Sect. \ref{conn}.  A number of other examples will be presented in \cite{Ho3}.

It deserves to point out that the map between continuous and discrete VOP allow in many cases easily to produce new examples from well known ones.

Although we use only fields of characteristic 0 many constructions could be performed in characteristic $p$.  In all such cases the polynomial system will be finite. It would be interesting to find what properties  these systems possess  and even their eventual applications. We recall that the Kravchuk polynomials, which form a finite orthogonal system, have applications to coding theory \cite{McWS}. 

The reader probably has noticed the absence in the above scheme of the Jacobi and Bessel polynomials as well as the polynomials of the Hahn family. I don't know if this is due to limitations of the method. However I have to point out that even the first case that could be considered to be an analogue of Jacobi or Bessel polynomials, i.e. when  the leading coefficient $a_3(x)$ of the operator $L = a_3(x) \partial_x^3 +\ldots$ is a polynomial of degree 3    (and $d=2$), is impossible, see \cite{HSTY}. This means that the Jacobi and Bessel polynomials are quite exceptional.

These investigations could be continued in several directions. First, one could further examine the properties of the VOP, introduced here.  We already pointed out that  the generalized Gould-Hopper polynomials  share other similarities with the classical orthogonal polynomials. The weights of these polynomial systems are subject to investigations in \cite{Ho4}. In the same paper we have constructed bi-orthogonal ensembles connected with  the generalized Gould-Hopper polynomials.

 Another study would be the  construction of new VOP via different versions of Darboux transformations in the spirit of \cite{GHH, GY, Il1,Il2}, etc. for the continuous families or the methods found by Dur\'an et al. \cite{Du, DdlI1, DdlI} in both  discrete and continuous case. 
 Notice that the methods of the present paper cannot produce these results - e.g. the Krall-Laguerre polynomials satisfy 3-term recurrence relation  with rational coefficients in the discrete variable $n$, (cf. \cite{GHH, Il2}),      while with the automorphisms of algebras method the finite-term recurrences  are  always with polynomial coefficients in $n$. More conceptual argument will be to compute the corresponding algebras of commuting differential operators, having the polynomial system as their eigenfunctions as in \cite{Il1, Il2}. In particular for the case Krall-Laguerre polynomials with measure  $e^{-x}  + M \delta(x)$ P. Iliev found that the algebra has 2 generators and more generally he obtained the number of generators for each of the measures. The same obviously holds for the discrete families. Contrary to this situation the polynomial systems obtained here are eigenfunctions of an algebra isomorphic to $\Cset[z]$.  This shows that both the methods of the present paper and the  cited above are relevant to the study of the GBP and none of them alone will give its complete solution.

It would be interesting to find applications to integrable systems along the lines of investigations in \cite{Cha}. However, it seems that  systems different from the  Novikov-Vesselov one would be the right object.  
In \cite{BHY2, HH, HHI} one can  find representations of the $W_{1+\infty}$-algebra and Virasoro algebra originating from certain solutions of the bispectral problem. It is  quite plausible that all solutions here have similar properties.  See   also \cite{Tak}  where string equations for certain solutions of the  2-d Toda hierarchy are derived.

The finite term recurrence relations generate   solutions for the bi-graded   Toda lattice hierarchy \cite{CDZ, Mi, MT,  BW}. It is worth understanding their nature.

\medskip

\noindent 
{\it Acknowledgements.} 
I am deeply grateful to Boris Shapiro for showing and discussing some examples of systems of polynomials studied here and in particular the examples from \cite{ST}. Without this probably the project would have never been even started.   The author is    grateful to the Mathematics Department   of Stockholm University for the hospitality in  April 2015. My thanks also go to the organizers and the participants of the International Conference Groups and Rings - Theory and Applications (GRiTA2016) and in particular  to Vesselin Drensky for the invitation to the conference,  where some of the results were presented. Recently I had the chance to present some of these results at the conference "On crossroads of analysis, algebra, and geometry or Boris' 60-th birthday, Stockholm. I thank the organizers for the invitation to present my results and the audience for the interest in my work. Also I am grateful to the participants of the Seminar "Dynamical systems and number theory" at Sofia University, where I presented a detailed account of the results.

The author is very grateful to the referees for careful reading of the manuscript and suggesting valuable improvements of the text.   

This research has been partially supported by the Grant No DN 02-5 of the Bulgarian Fund "Scientific Research".

\part{Theory}\label{I}


\section{Elements of bispectral theory} \label{BT}

Let $\R_1$ and $\R_2$ be associative algebras over a field $k$ and let $\M$ be a
left module over both of them. Let $\F_1$ and $\F_2$ be fields such that
$k\subset\F_1\subset\R_1$ and $k\subset\F_2\subset\R_2$.
\bde{1.1}
We call an element $L\in\R_1$ {\em bispectral\/} iff there exist $\psi\in\M$,
$\Lambda\in\R_2$, $f\in\F_2$, $\theta\in\F_1$ such that
\beqa
L\psi=f\psi,
\label{1.1}\hfill\\
\Lambda\psi=\theta\psi.
\label{1.2}\hfill
\eeqa
\qed
\ede
In order to have a non-trivial problem we assume that

\begin{quotation}
$(**)$ \textit{there are no nonzero elements $L$ and $\La$
that satisfy one of the above conditions with $f \equiv 0$ or $\theta \equiv 0$.}
\end{quotation}

Let us fix $\psi\in\M$ satisfying \eqref{1.1} and \eqref{1.2}. We are
interested in the equation (cf.\ \cite{BHY})
\beq
P\psi=Q\psi
\label{1.4}
\eeq
for $P\in\R_1$, $Q\in\R_2$. We put
\begin{eqnarray*}
&&\B_1 = \{ P\in\R_1 \; | \; \exists Q\in\R_2 {\textrm{ for which \eqref{1.4}
		is satisfied}} \},
\hfill\\
&&\B_2 = \{ Q\in\R_2 \; | \; \exists P\in\R_1 {\textrm{ for which \eqref{1.4}
		is satisfied}} \}.
\hfill
\end{eqnarray*}
We shall assume that the actions of $\R_1$ and $\R_2$ on $\M$ commute and that (**)
  holds. Then $\B_1$ and $\B_2$ are associative algebras over
$k$ without zero divisors. Obviously \eqref{1.4} defines an anti-isomorphism
\[
b\colon \B_1 \to \B_2, \quad b(P)=Q.
\]
Introduce also the subalgebras
\begin{eqnarray*}
\K_1 = \B_1\cap \F_1, \quad \K_2 = \B_2 \cap \F_2,
\hfill\\
\A_1 = b^{-1}(\K_2), \quad \A_2 = b(\K_1).
\hfill
\end{eqnarray*}
Then $\A_1, \A_2 $   are a commutative algebras isomorphic to $\K_2, \K_1$ (the {\em
	spectral algebras\/}).  
\bde{1.2}
We call the triple $(\B_1, \B_2, b)$ a {\em{bispectral
		triple}} iff both $\K_1$ and $\K_2$ contain non-zero
elements.
\qed
\ede
Obviously if $(\B_1, \B_2, b)$ is a bispectral  triple then any
element $L \in \A_1$ is bispectral.

Examples of such situations with differential operators are given in the Introduction. Let us consider other examples.

\bex{dOP}
In this example  we consider functions $f(x)$ which are invariant under the shift operator $D$, acting as $Df(x)  = f(x+1)$ and its inverse $D^{-1}$, $D^{-1}f(x)  = f(x-1)$. The operators $D^{\pm}$     and the operator $\hat{x}$ of multiplication by $x$ form the  algebra $\R_1$ in this case. We are going to introduce $\R_2$ exactly in the same way considering functions $h(n)$ and the corresponding shift operators $T^{\pm}$ as well as the operator $\hat{n}$ of multiplication by $n$. Recall the definition of the falling factorial $(x)_n$ 

\[
(x)_n = x(x-1)\ldots (x - n +1),\; n\in \Nset, \; (x)_0 = 1.
\]

\bre{ff}
This notation is not the usual one in the research in special functions. Usually it denotes the rising factorial. In this paper and in \cite{Ho3} we choose to use the above notation for the notion of the falling factorials. While we need both the rising and the falling factorials we use more often the latter.

\ere 
Let us put 

\[
\psi(x,n) = (x)_n.
\]
It would be convenient to introduce also the operators $\Delta = D- 1$ and $\nabla = D^{-1} - 1$ and $H = -\hat{x}\nabla$.
We are going to denote  by $\B_1$ the subalgebra of $\R_1$ generated by the   operators   $\Delta$, $H$ and the operator   $\hat{x}$. 
Then we have

\[
\begin{cases}
\Delta \psi(x,n) = nT^{-1}\psi(x,n)\\
x \psi(x,n) =     (T  + n)\psi(x,n)\\
H\psi(x,n) = n\psi(x,n).
\end{cases}
\]
This defines an anti-involution $b$, defined  by $\psi(x,n)$:

\[
\begin{cases}
b(\Delta)   = nT^{-1}\\
b(x) =     (T  + n)\\
b(H)= n.
\end{cases}
\]
The subalgebra $\B_2$ is  generated by the image of $\B_1$ under the action of the bispectral involution $b$  defined by $\psi(x,n)$

\eex
\qed

Of course the  most important examples are the classical orthogonal polynomials of Hermite, Laguerre, Bessel and Jacobi, see e.g. \cite{KLS, NIST}.



In the next proposition we will discuss how to construct new bispectral operators from already  known ones. 
Our basic  tool for constructing new bispectral operators  will be the following simple observation made in \cite{BHY}.
\bpr{2.1} Let $(\B_1, \B_2, b)$ be a bispectral triple and  let    $L \in \B_1$ be locally nilpotent. Suppose that,  the element  $e^L\psi$ is well defined, $e^L\psi \in \M$. Define a new map $\sigma = e^{\ad_L} :\B_1 \rightarrow \B_1$ via the new  element $\psi'  := e^L\psi$. 

Then $b^{'}: \B_1\to \B_2$ is surjective, and $(\B_1, \B_2, b^{'}  )$ is a bispectral triple. (Here {\em locally nilpotent\/} means that for each $A\in\B_1$ there exists an
integer $n$ such that $(\ad_L)^n A = 0$.)
\epr

\smallskip\noindent

\section{Algebraic construction} \label{AC}

\subsection{The Weyl Algebra and some of its subalgebras} Below we repeat the definition of 
 the Weyl algebra $W_1$ in the abstract form in which  we are going to use it. $W_1$ is an algebra, generated over a field $\F$ by two elements $Y, Z$, subject to the relation $[Z, Y] \equiv ZY-YZ = 1$.    Denote by $H$ the element $YZ$ and let $R(H)$ be a polynomial of degree $d$ in   $H$. Consider also the element $G$ defined as 
\beq
G= R(H)Z \in W_1.
  \label{G-def}
\eeq
We are  going to define  a subalgebra $\B^R_1$ of   $W_1$,   generated by the elements $Y, H, G$.

In the next lemma we list    the relations between these elements. We have

\ble{c-r}
The elements $ Y, H, G$ commute as follows:

\[
[H, Y]= Y,  \;\; [H, G] =   - G, \;\; [G, Y]=  R(H)(H+ 1) -  R(H-1)H.
\]
\ele
\proof Follows from the defining relations of $W_1$ by direct computations. 

\qed

\bigskip

 For later use  we put $S(H) := R(H-1)H$. Then $[G, Y]=  \Delta S(H) := S(H+1)- S(H)$.

Next we define an automorphism $\sigma$ of $\B^R_1$, acting on elements $A\in \B^R_1$  as 

\[
\sigma(A) = e^{\ad_G} (A)  =  \sum_{j=0}^{\infty} \frac{\ad^j_G A}{j!},
\]
see \cite{Dix}. This definition makes sense only if the above series is finite for any   $A\in \B^R_1$  (the number of terms depends on $A$).
Let us compute the action of $\sigma$, which will also prove the correctness of $\sigma$. 

\ble{au}
The automorphism $\sigma$ acts on the generators of $\B_1$ as follows:

\begin{equation*}
\begin{cases}
\sigma(G) = G\\
\sigma(H) = H + G\\
\sigma(Y) =  Y      +  \sum\limits_{j=1}^{d+1}\frac{[\Delta^j S(H)]G^{j-1}}{j!}. \\
\end{cases}   
\end{equation*}
\ele
\proof

The  first two formulas are obvious.  Let us prove the last one. We start with the commutator

\[
\ad_G(Y) = [G, Y]  = \Delta S(H). 
\]
Let us derive  the general formula for $\ad_G(f(H))$ for any polynomial $f(H)$. We have 

\[
[G, f(H)] = f(H+1)G - f(H)G = \Delta f(H)G.
\]
From this we find

\[
\ad_G^2(Y)  =   [G, \Delta S(H)]= \Delta^2 S(H)G.
\]
By induction we have

\[
\ad_G^j( Y)= [\Delta^j S(H)]G^{j-1}.
\]
On the other hand $\Delta^j S(H)$ is a polynomial of degree $d+1 -j$. Hence $\ad_G^{d+2}(Y) = 0$. This gives the last formula. 

\qed

We can define more  general automorphisms by making use of any polynomial $q(G)$ in $G$ without a constant term   and put 

\[
\sigma_q = e^{\ad_{q(G)}}.
\]

\bre{constant}
Without the restriction the class of new VOP will be the same. The point is that a nonzero constant term $q_0$ in $q(G)$ will simply multiply the polynomials by the constant factor $e^{q_0}$. The 
restriction is made only to simplify the formulas and the arguments. The same restriction will be imposed through the entire paper.
\ere 

In what follows we will drop the index $q$ from $\sigma_q$ as there is no danger of confusion.

\ble{gen2}
$\ad_{q(G)}$ acts locally nilpotently on $\B_1^R$, i.e. the series defining $\sigma_q(A)$ is finite for any $A$.
\ele
\proof 
We are going to prove the statement for the generators   $Y, H, G$   of $\B_1^R$.
First,     $\sigma(G) = G$.

Next we consider $\sigma(H)$. We have 

\beq \label{G-H}
[G^k, H] = k G^k.
\eeq
Hence

\[
\ad_{q(G)} (H)    = q'(G)G  \; \; \text{and} \;\; \ad^2_{q(G)}(H)= 0.
\]
This gives

\[
\sigma(H) = H + q'(G)G.
\]

Finally let us  compute

\[
\sigma(Y) =  Y + \sum_{j=1}^{\infty} \frac{\ad^j_{q(G)}(Y) }{j!}.
\]  
We start with  $[q(G), Y]$. From the identity

\[
[G, Y] = \Delta S(H)
\]
we find for any $k$

\[
[G^k, Y] = \sum_{j=0}^{k-1} G^j\Delta S(H)G^{k-1 -j}.
\]
In what follows we do not need the exact coefficients at the powers of $G$.  Using \eqref{G-H} we can write the last formula as 

\[
[G^k, Y] = F(H) G^{k-1}, \; \deg (F) = d.
\]
Let us denote the degree  of $q(G)$ by $l$. Then we see that $\ad_{q(G)}$  has the form

\[
\ad_{q(G)}(Y) =  \sum_{j <l} Q_j (H)G^j, \;\; \deg Q_j  = d.
\]
Each application of $\ad_{q(G)}$ will reduce the degree of  $Q_j$ by 1.  From \eqref{G-H} with $k=l$ we see that it raises the degree of $G$ by $l$. Hence a  formula of the form

\[
\ad^r_{q(G)}(Y) =  \sum_{j < lr }  Q_{d-r +1}(H) G^j 
\]
holds with $\deg Q_{d-r +1}(H) = d-r +1$. This yields

\[
\ad^{d+2}_{q(G)}(Y) = 0.
\] 

\qed

From the proof it follows that 
\bco{rec-q}
The automorphism $\sigma$ acts on $Y$ as

\[
\sigma Y = Y + \sum\limits_{j=0}^{ld+l-1}\ga_j(H)G^j.
\]
\eco

\qed

The automorphisms induced by $q(G)$  will    also be used in our studies of VOP. 
 
\subsection{Bispectral anti-involution}
Let $\R_2$ be an  algebra with unit $1$, generated by the generators $T, T^{-1}, \hat{n}$ and the relations

\[
T \cdot T^{-1} = T^{-1} \cdot T = 1,\;    [T, \hat{n}]  =  T, \;\; [T^{-1}, \hat{n}] = - T^{-1}.
\]
It is easy to see that $[T, \hat{n} T^{-1}] = 1$. Hence the elements $T, \hat{n} T^{-1}$ span a copy of the Weyl algebra.
Let us define an anti-homomorphism $b$, i.e. a map $b\colon \B^R_1 \to \R_2$ satisfying  $b(m_1.m_2) = b(m_2).b(m_1)$), for each $m_1,\; m_2 \in \B_1$  by:

\[
\begin{cases}
b(Y)= T\\
b(H) = \hat{n}\\
b(G) = \hat{n} T^{-1}R(\hat{n}).
\end{cases}
\]
 Let us define the algebra $\B_2^R = b(\B_1^R)$. Then  $b\colon \B^R_1 \to \B_2^R$ is an anti-isomorphism and in particular $b^{-1} \colon      \B^R_2 \to \R_1$ is well defined.

\bre{pol-abs}
In section \ref{BT} we first defined a bi-module $\M$ and using an element $\psi \in \M$ we defined the anti-isomorphism  $b$. As we see there is no problem to perform everything without $\M$. Nevertheless, having in mind our goal of constructing VOP, we can point out a suitable bi-module and an element $\psi$ in it.  Let us take $\M = \F[X]$,  where   $X$ is transcendental over $\F$ and put for example $\psi = X^n, \; n = 0, 1, \ldots$. Then we can define actions on $\psi$ by

\[
\begin{cases}
Z(X^n)= nX^{n-1}\\
Y(X^n)= X^{n+1}\\
T(X^n)= X^{n+1}\\
\n T^{-1}(X^n)= nX^{n-1}.\\
\end{cases}
\]

\ere
\qed

\subsection{New anti-isomorphism}

Using the automorphism $\sigma$ we can define a new anti-isomorphism $b'$ by the formula $b' = b\circ \sigma^{-1}$. In what follows we will compute explicitly $b'$. The more general $b'_q:= b\circ \sigma_q^{-1}$ can also be computed explicitly. However the formulas are too involved and in any specific situation they are easier to compute directly. Below having in mind applications to VOP we will use the generator $L = \sigma(H)$ instead of $H$.

\bth{th:new}
	The new bispectral involution $b^{'}$ on $\B_1$ satisfies: 
	
	\begin{equation}
	\begin{cases}
		b^{'}(G)  = \hat{n}R(\hat{n}-1) T^{-1}\\
	 b'(Y)  = T+  \sum\limits_{j=1}^{d+1}\frac{(-1)^j[R(\hat{n}-1)nT^{-1}]^{j-1} [\Delta^j S(\hat{n})]}{j!} \\
		b^{'}(L)  =\hat{n}.
	\end{cases}  \label{eq:nbi}
	\end{equation}
\ethe 

\proof
We have 
$b^{'}(G)  = b \circ\sigma^{-1} (G) = b (G) = \hat{n}R(\hat{n}-1)T^{-1}.$ 
In the same way we compute
\begin{eqnarray*}
	b^{'}(Y)  = b \circ\sigma^{-1} (Y) = b(Y) + \sum\limits_{j= 1}^{d+1}\frac{(-1)^jb([\Delta^j S(H)]G^{j-1})}{j!}    \\
	=  T+ \sum_{j=1}^{d+1}\frac{(-1)^j(\hat{n}R(\hat{n}-1)T^{-1})^{j-1}[\Delta^j S(\hat{n})]}{j!}.  
\end{eqnarray*}
Finally for $L$ we find
$$
b^{'}(L)  = b\circ \sigma^{-1}(L)  = b\circ \sigma^{-1}    \circ\sigma(H) = b(H) = \hat{n}.
$$
\qed

Following  \cite{Ho1, Ho2} we introduce  more general automorphisms. Let $q(Y)$ be a polynomial without a constant term. We can define an automorphism $\sigma$ as

\[
\sigma_q = e^{\ad_{q(G)}}\colon \B_1^G \to \B_1^G.
\]
From it we define the anti-automorphism $b^q = b\circ\sigma_q^{-1}$. Then the following theorem holds

\bth{th:n-gen}
The new bispectral involution $b^{'}$ on $\B_1$ satisfies: 

\begin{equation*}
\begin{cases}
b^{'}(G)  = \hat{n}R(\hat{n}-1)T^{-1} \\
b^{'}(L)  = \hat{n}
\end{cases}  
\end{equation*}
and an equation of the form

\[
b'(Y)  = T+  \sum\limits_{j=1}^{ld-1}\gamma_j(H)(\hat{n}T^{-1})^j, 
\]
 where coefficients $\gamma_j$ are polynomials in $H$.
\ethe 

\proof  Essentially repeats  the proof of   \thref{th:new}.
 
\qed

\part{Applications} \label{II}

\section{Continuous VOP}\label{con}

In this section we want to use the abstract algebraic construction to produce VOP, which are eigenfunctions of ordinary differential operators. With this purpose in mind     we realize the Weyl algebra $W_1$  as an algebra of differential operators in one variable $x$ with  polynomial coefficients. We consider that the Weyl algebra is represented in $\F[x]$ by realizing    $Y$ as the operator of multiplication  by $x$ and  the  operator     $Z$ as the differentiation  $\partial_x$ with respect to $x$.

  The operators  $T^{\pm}$ will be acting on functions $f(n)$ by  $T^{\pm}f(n) = f(n \pm 1)$ and the operator $\hat{n}$ as   multiplication by $n$, i.e. $\hat{n} f(n) = nf(n)$. We can apply  the algebraic constructions   from Part \ref{I}  to this case.  In particular we define  $H= x\partial_x$. We fix a polynomial $R(X)$ and define the element $G = R(H)\partial_x$. Then our  algebra $\B_1^R$ will be generated by $G, H, \hat{x}$.

Our initial polynomial system will be     $\psi(x, n) = x^n$. Then we have

 \[
 \begin{cases}
 G \psi(x,n) = \hat{n}R(\hat{n}-1) T^{-1} \psi(x,n)\\
 H \psi(x,n) = \hat{n} \psi(x,n)\\
 x \psi(x,n) =  T \psi(x,n).
 \end{cases}
 \]
   Next we define the corresponding automorphism 
 
 \[
 \sigma = e^{\ad_G}.
 \]
With its help we define the vector orthogonal polynomials. This will be done by putting

\[
P_n^R(x) = e^{G}\psi(x,n)= \sum_{j=0}^{\infty}\frac{G^j\psi(x,n)}{j!}.
\]
Notice that the operator $G$ reduces the degree of any polynomial by 1.  This yields that  the above series is finite and defines a polynomial of degree $n$.

Let us introduce the operator  $L = \sigma(H)$. In the present situation it reads

\[
L = x\partial_x + R(x\partial_x)\partial_x.
\]
Finally we define the differential operator $M= \sigma(x)$.

\bth{2nd}
	The polynomials $P_n^R$ have the following properties:

	(i) They are eigenfunctions of the differential operator  $L$
	 with eigenvalues $\lambda(n) = n$.
	
	(ii) They satisfy the recurrence relation
	
	\begin{eqnarray}\label{new-rec}
		x P_n^R(x)   =       b'(x) =   P_{n+1}^R(x)  +   \sum\limits_{j=1}^{d+1}\frac{(-1)^j[nR(n-1)T^{-1}]^{j-1} \Delta^jS(n)}{j!} P_n^R(x).
	\end{eqnarray}
	
	(iii) They have ladder operators as follows: the operator $G$ is a lowering operator

	\begin{equation} \label{new-low}
	G P^R_n(x)  = nR(n-1)P^R_{n-1}(x) 
		\end{equation}
	while the operator $M$ is a raising operator:
	
	\begin{equation} \label{new-rai}
	M P_n^R(x) = P_{n+1}^R(x).
	\end{equation}

(iv) Rodrigues-type formula holds:

	\[
  P_n^R(x)  =	M^n 1.
	\]

\ethe 

\proof The proof follows quite directly from \thref{th:new}  Indeed, statement (i) that   polynomials $P_n^R(x)$ are eigenfunctions of the operator
	\[
	L =   R(H)\partial_x+ x\partial_x.  
	\]
	follows from the computation of the new bispectral involution $b'$ in  \thref{th:new}.  More exactly, the third equation of \eqref{eq:nbi}  gives
	\[
 (R(H)\partial_x+ x\partial_x) P_n^R(x)  = n P_n^R(x).
	\] 
	To prove recurrence (ii), observe that the second equation in \eqref{eq:nbi} provides 
	\[
	x P_n^R(x)   =       b'(x) =   P_{n+1}^R(x)  +   \sum\limits_{j=1}^{d+1}\frac{(-1)^j[nR(n-1)T^{-1}]^{j-1} \Delta^jS(n)}{j!} P_n^R(x).
	\]
	This gives the desired formula \eqref{new-rec}.

	Also the third equation in \eqref{eq:nbi}  gives the lowering operator \eqref{new-low}
	\[
	R(H)\partial_x P_n^R(x) = nR(n-1)P_{n-1}^R(x).
	\]
Finally we have

\[
   M P_n^R(x)  = e^Gxe^{-G}  e^G x^n =  e^Gx x^n = P_{n+1}^R(x),
\]	
which proves \eqref{new-rai}.
The Rodrigues-like formula follows from the last one.

\qed

Next we  consider a   polynomial $q(X)$ and  we assume that
$q(X)$ has no constant term.  Above we formulated the special case of $q(G) = G$ as the formulas there are simple.

 Define an automorphism of the   algebra $\B_1$

\[
\sigma(A) = e^{\ad_{q(G)}}(A),\;\; A \in \B_1.
\]
Also we define 

\[
P_n^q(x)  = e^{q(G)}x^n   = \sum\limits_{j=0}^{\infty}\frac{(q(G))^j x^n}{j!}.
\]
It is clear that the latter series contains only a finite number of terms as the operator     $q(G)$    reduces the degrees of the polynomials. (Here we use that $q$ has no constant term.)   Hence $\psi(x,n)$ is a polynomial in $x$. 
Define  the operator  $L: = \sigma(H)$. One can easily see that  

\[
L = x\partial_x +  q'(G)G. 
\]
%

We formulate the  main results about the polynomial system $P_n^q(x)$ in the next theorem. 

\bth{any}
The polynomials $P_n^q(x)$ have the following properties:

(i) They are eigenfunctions of the differential operator 

\begin{equation*}      
	L :=  q^{'}(G)G  + x\partial_x.  
\end{equation*}

(ii) They satisfy    a    recurrence    relation  of the form

\begin{equation}     \label{new-rec3}
xP^q_n(x)  =P^q_{n+1} + \sum_{j=0}^{ld +l-1}\gamma_j(n)P^q_{n-j}. 
\end{equation}

(iii) They have ladder operators as follows. The operator  $G$ is a lowering operator:

\begin{equation*}         
	G P^q_n(x)  = nR(n-1) P^q_{n-1}.   
\end{equation*}
while the    operator $M = \sigma_q(x)$ is a raising operator:

\[
MP^q_n(x) = P^q_{n+1}(x).
\]
\ethe

(iv) Rodrigues-like formula holds:

\[
P_n^q(x)  = M^n. 1.
\]
We skip the proof as it   essentially repeats the proof of \thref{2nd}.

\bre{RR}
Notice that the order           of the differential operator $q(G)$ is $l(d+1)$. The recurrence relation involves  $l(d+1) +1$   terms. 
This fact  is similar to  the case  of the classical orthogonal     polynomials, where the differential operator is of order 2, while the recurrence relation is  3-term. The remark answers a question by B. Shapiro (private communication). 
\ere
 
\section{Discrete VOP} \label{dis}

Naturally here we use a realization of the Weyl algebra in terms of difference operators 
  acting on   the space    of polynomials $\F[x]$. Then we  define the operators   $D^{\pm}$ acting on $f \in \F[x]$ by shift of the argument   $D^{\pm}f(x) = f(x \pm 1)$. The operator $\hat{x}$ acts on $f(x)$ as multiplication by $x$. We also need      $\Delta = D -1$,  $\nabla = D^{-1} - 1$ and $H = - \hat{x}\nabla$. Finally we put $g = \hat{x}  - H$. The proofs here are almost the same as in the  continuous case   and we only point the differences.

  We formulate the needed relations in the following lemma.  
  
  \ble{dis-W}
  The operators $g, \Delta$ span an algebra isomorphic to the Weyl algebra. The relations, between the operators defined above, are as follows:
  
  \[
  \begin{cases}
  [\Delta, g] = 1\\
  [\Delta, H] = \Delta\\
  [H, g] = g.
  \end{cases}
  \]
  \ele
  \proof Simple computations.
  
  \qed

  Let  $\psi(x,n):=  (x)_n$. The next lemma shows that $\psi(x,n)$ realizes the bispectral involution.
  
  \ble{dis-bi} 
  The following identities hold:
  \[
  \begin{cases}
  H\psi(x,n) = n\psi(x,n)\\
  g\psi(x,n) = T \psi(x,n)\\
  \Delta\psi(x,n) = n T^{-1}\psi(x,n).\\
  \end{cases}
  \]
  \ele
  
  \proof Simple computations.
  
  \qed
  
  We fix a polynomial $R(H)$ and define the element $G = R(H)\Delta$. Then our  algebra $\B_1^R$ will be generated by $G, H, \hat{x}$. Instead of $\hat{x}$ it will be more convenient to use $g$.      Then we have

  \[
  \begin{cases}
  G \psi(x,n) = nR(n-1) T^{-1} \psi(x,n)\\
  H \psi(x,n) = n \psi(x,n)\\
  g \psi(x,n) =  T \psi(x,n).
  \end{cases}
  \]

We see that all conditions to apply the algebraic  construction are met.  We can define an automorphism $\sigma\colon \B_1 \to \B_1$  as  

\[
\sigma = e^{\ad_G}.
\]

\ble{dis-gen}
The automorphism $\sigma$ maps the generators of $\B_1^R$ as follows

\[
\begin{cases}
\sigma(g) = g     +  \sum_{j=0}^{d}\frac{\Delta^{j+1} S(H)G^j}{j!} \\
\sigma(H) = H + G \\
\sigma(G) = G.\\
\end{cases}
\]
\ele

\proof It is a special case of \leref{au}.

\qed

Notice that  we  need the image of $x$ instead of $g$. It is easily obtained by 

\[
\sigma(x) = \sigma(g) + \sigma(H) = x +  \sum_{j=0}^{d}\frac{\Delta^{j+1} S(H)G^j}{j!}
 + G. 
\]

Let us introduce the polynomials $P_n^R(x)$ by

\[
P_n^R(x) = e^{G}\psi(x,n) = \sum_{j=0}^{\infty}  \frac{G^j \psi(x,n)}{j!}.
\]
Again using the argument  that the operator $G$ reduces the degree of any polynomial by 1 we find  that the series is finite and $P_n^R(x)$ is a polynomial of degree $n$.      Also we define the operator

\[
L= \sigma(H) = H + G.
\]

Now we are in a position to apply the   algebraic construction of the new bispectral involution.
\ble{dis-new-b}
The new bispectral involution $b' = b\circ \sigma^{-1}$ acts on the generators of
$\B^R_1$ as follows

\[
\begin{cases}
b'(x)  = T +n  +    nR(n-1)T^{-1} +  \sum\limits_{j=1}^{d+1}\frac{(-1)^{j+1}   (nR(n-1)T^{-1})^j \Delta_n^{j+1}S(n) }{j!} \\
b'(L)  =  n   \\
b'(G)  =   nT^{-1}R(n). \\
\end{cases}
\]
\ele
\proof   \leref{dis-new-b} is just  \thref{th:new}   in the specific situation.

\qed

Exactly as in the continuous case we can define automorphisms using polynomials $q(G)$  of $G$ without a  constant term. Then we can define the automorphism
\[
\sigma_q = e^{\ad_{q(G)}}.
\]
The corresponding polynomials that we need are

\[
P^q_n(x)   = e^{q(G)}(x)_n.
\]
Let us define the difference operator

\[
L = x\nabla + q'(G)G.
\]
In the next theorem we formulate the  main results about the polynomial system $P_n^q(x)$. The special case of $q(G) = G$ will be only commented.

 \bth{any-dis}
 The polynomials $P_n^q(x)$ have the following properties:

 (i) They are eigenfunctions of the difference operator 
 
 \[
 L :=  q'(G)G  - x\nabla.
 \]
 
 (ii) They satisfy    a    recurrence    relation  of the form 
 
 \[
 xP^q_n(x)  =P^q_{n+1} + \sum_{j= 0}^{m}\gamma_j(n)P^q_{n-j}. 
 \]
 where $m= ld+l-1$ when $d>0$ and $m = l$  when $d=0$.
 
 (iii) 	  They possess ladder operators as follows. The operator $G$ is lowering operator:
 
 \[ 
 G P^q_n(x)  = nR(n-1) P^q_{n-1}; 
 \]  
 the operator 
 $M = \sigma(x)$ is a raising operator:
 
 \[
 M P_n^q(x) = P_{n+1}^q(x).
 \]
 
 (iv) Rodrigues type formula holds:
 
 \[
 P_n^q(x) = M^n.1.
 \]

 \ethe
 
 \proof The proof is essentially the same as the proof of \thref{2nd}. The only difference comes from the recurrence relation. The point is that $b'(x) = b'(g) + b'(H)$. While the term $b'(g)$ is the same as   $b'(x)$  in \thref{2nd}, the other term $b'(H)$ contributes to $b'(x)$ with $q'(G)G$. Its  degree in  $G$  is always $l$.  The degree in $G$ of  $\ad^{d+1}(g)$ is always $ld +l -1$. But for $d=0$ this gives $l-1$, which is less than $l$. For $d>0$  we have $ld +l -1 \geq l$. 
See Part \ref{III} for examples.
 
 \qed
 
 We only mention the explicit form of the recurrence relation in the case $q(G) = G$:

\[
\begin{split}
xP_n^R(x) &= P_{n+1}^R(x) +nP_n^R(x)  + nR(n-1)P_{n-1}^R(x)       +\\
& +  \sum_{j=1}^{d+1}\frac{(-1)^j\prod_{k=1}^{j-1}(n-k+1)R(n - k)[\Delta^j S(n-j)]}{j!}P^R_{n-j+1}(x).
\end{split}
\]

\section{Discrete vs. continous VOP} \label{DdsC}

The common origin of the VOP in the continuous and the discrete cases suggests that there is quite direct   connection between them pairwise. Let us fix the polynomial $R(H)$. For the purposes of this section we denote the algebras of differential operators corresponding to $R$ by $\B^{c}_1$, while the algebra of difference operators will be $\B^{d}_1$, supressing the dependence on $R$. We also will attach the corresponding indexes $c$ or  $d$ to the corresponding  operators, as follows $A \in \B^{c}_1$ will be $A^c$, $A \in \B^{d}_1$ will be $A^d$. 

First we point the connection between the basic systems. From the definition of the gamma function

\[
\Gamma(-x + n) = \int_{0}^{\infty} t^{-x -1 +n}e^{-t}dt 
\]
and using that

\[
(x)_n =(-1) ^n \frac{\Ga(-x + n)}{\Ga(-x)}
\]
we obtain

\[
(-1)^n(x)_n = \frac{1}{\Gamma(-x)}\int_{0}^{\infty} t^{-x -1 +n} e^{-t}dt.
\]
This prompts the definition of Mellin type transform of the function $f(t)$

\[
\mathcal{M}^*(f)(x) = \frac{1}{\Gamma(-x)}\int_{0}^{\infty} f(t) t^{-x-1} e^{-t}dt.
\]
The main feature of  the transform $\mathcal{M}^*$ pertinent here is that it transforms the system $t^n$ into the system $(x)_n$:

\[
\mathcal{M}^*(t^n) = (x)_n.
\]
\bre{Mel}
In fact if we consider Mellin transform of $e^{-t}f(-t)$ then 

\[
\mathcal{M}((-t)^ne^{-t})(-x) = \Gamma(-x)(x)_n.
\]
We are going to use the slight modification $\mathcal{M}^*$ defined above to underline the isomorphism between two copies of the vector space of polynomials $\Cset[x]$.

\qed

Notice however that we don't need to go to integration and to the use of transcendental functions. We use Mellin's transform only to point the origin of the map. Of course one can   define the map in purely algebraic terms over any field via the Stirling numbers $s(n,k)$, which   sometimes  are  introduced  also  by the isomorphism of $\F[x]$
given by
\[
(x)_n = \sum\limits_{k=0}^{n}s(n,k)x^k.
\]
This shows that we again can work with any field.
 \ere
 
 \bpr{tran} Let us fix the polynomial $R(H)$
 The   Mellin type transform $\mathcal{M}^*$ maps the continuous VOP $\{P_n^{c}\}$  into the discrete ones with th same $R$   $\{P_n^{d}\}$. Moreover
 
 \[
 \mathcal{M}^*(P_n^{c}(t) ) = P_n^{d}(x).
 \]
 
 \epr
 \proof 
 Let $G$ be the operator \eqref{G-def}. We apply its version in the usual basis $\{t^n\}$.   We have 
 
 \[
 (G^c)^j t^n = (n)_j\prod_{k=1}^d(n-1+\alpha_k)_jt^{n- j}.
 \]
 Obviously  
 
 \[
 \mathcal{M}^*((G^c)^j t^n)= (n)_j\prod_{k=1}^d(n-1+\alpha_k)_j\mathcal{M}^*(t^{n-j})        = (G^d)^j (x)_n. 
 \]
 Then applying  $\mathcal{M}^*$  to the definition of the polynomials $\{P_n^{c}\}$
 
 \[
 P_n^{c}(t) = \sum\limits_{j=0}^{\infty}\frac{(G^c)^j t^n}{j!}
 \]
 and plugging in it the expressions for $(G^c)^j t^n$ we obtain
 
 \[
  \mathcal{M}^*(P_n^{c}(t)) = \sum\limits_{j=0}^{\infty}\frac{(G^d)^j (x)_n}{j!} = P_n^{d}(x).
 \]
    
 \qed

 Of course we can extend the map $\mathcal{M}^*$ to include the cases of $q(G)$. This is straightforward. We only mention the special case of 
 the polynomials that originate from $q(\partial_x)$, which are usually called 
  Appell    polynomials (see, e.g. \cite{Dou}). It is natural to call the polynomials defined via $q(\Delta)$ Charlier-Appell polynomials, cf. \cite{Ho2} and  \exref{ex1-d}. In the case when $q(G) = a\partial_x^l, \; a \in \Cset$,   the polynomials are the well known Gould-Hopper polynomials \cite{GH}. 
 
 The next corollary is quite straightforward but we formulate it explicitly because of   its importance.
  
 \bco{pairs}
 The following systems are mapped to each other:
 
 (i)   The $d$-orthogonal Appell polynomials are mapped to $d+1$-orthogonal Charlier-Appell polynomials.  In particular the image of the classical  Hermite polynomials are 2-orthogonal polynomials. The Charlier polynomials are obtained from the first member (trivial) of the Appell family, i.e. with   $L = x\partial_x + a\partial_x$. 
 
 (ii)   The classical Laguerre polynomials are mapped onto the Meixner polynomials.
 \eco
 
 \proof   We have left the arguments for    Part \ref{III}.
 
\part{Examples} \label{III}
  
 \section{Elementary examples} \label{ee}
 In the examples below we construct the differential operator $L$ that has the VOP   $P_n^q(x) = e^{q(G)}x^n$ (or the corresponding expression $P_n^q(x) = e^{q(G)}(x)_n$ in the discrete case) as eigenfunctions. In some examples  we also write explicitly the finite-term recurrence relation.

 \subsection{Continous VOP}

 \bex{Seg} Set $q(\partial_x) =     -\partial_x^l/l $. Then we obtain
 \[
 L =     x\partial_x   -   \partial_x^l. 
 \]
 Observe that, for $l=2,$   the operator $ -L$ is  the standard Hermite operator.   The recurrence relation for $P_n(x)$ is given by \eqref{new-rec}. In our case, it reads:
 \[
 xP_n(x)  =  P_{n+1}(x)  + n(n-1)\ldots (n - l +2) P_{n-l+1}(x),
 \]
 which agrees with the results in \cite{ST}.   If $q(\partial_x) =     \rho \partial_x^l$ these are the well known Gould-Hopper polynomials, see  \cite{GH}. They are discussed in more detail in \cite{Ho3}.
 
    In general all the examples of with $G=  \partial_x$ and arbitrary $q$    have been studied thoroughly starting with the classical paper by P. Appell \cite{App} and bear his name. In the context of vector orthogonal polynomials the results have been found in \cite{Dou}.
 
 \eex

 \bex{ex2-1}
 Here we take $G= -(x\partial_x +\al + 1)\partial_x$. Let us take the polynomial $q(G) = G$.
 In this case we have
 
 \[
 L_1 = H + q'(G)G = -x\partial_x^2 -(\al +1-   x)\partial_x
 \]
 which   is the generalized Laguerre operator whose eigenfunctions are the Laguerre polynomials $L^{(\al)}_n(x)$.  See also \exref{L-M}.
  \eex
 
 \bex{ex2-2}
 The simplest new example is given by the polynomial 
 $q(G) = G^2/2$,    $G = x\partial_x^2 +\beta \partial_x$.  
 In this case we have

 \begin{eqnarray*}
 	L & =&  H + q'(G)G =  (x\partial_x^2 +\beta \partial_x)^2   + x\partial_x  \\
 	&= &x^2\partial_x^4 + 2x(1 + \beta)\partial_x^3  +  (\beta + \beta^2)\partial_x^2 + x\partial_x.
 \end{eqnarray*}
 Using the commutation relations between $\hat{x}, G$ and $H$ we can compute that 
 
 \[
 \sigma(x) = x + H G + GH +  G^3,
 \]
 which gives 
 
 \[
 b'(\hat{x}) = T - (n)_2n(n +\beta)(2n-1)T^{-1}  +(n)_3(n+\beta  -1)_3T^{-3}.
 \]

 The 5-terms recurrence relation reads
 
 \[
 xP_n(x) = P_{n+1}(x) -  n(n +\beta-1)(2n-1) P_{n-1}(x) + (n)_3(n -1 + \beta)_3P_{n-3}(x).
 \]
 Here are the first few polynomials:
 
 \[
 P_0(x) = 1, P_1(x) = x, P_2(x) = x^2 +2\beta(1 +\beta), P_3(x) = x^3 + 2.3.(1 +\beta)(2 +\beta)x,\ldots 
 \]
  A simple inspection shows   that the polynomials with even indexes are even and the ones with odd indexes are odd. Put $x^2 = u$. Let us introduce the following new polynomial systems:
 
 \[
 Y_n (u)= P_{2n}(\sqrt{u}), \; Z_n = u^{-1/2} P_{2n+1}(\sqrt{u}).
 \]
 From the 5-term recursion relation we derive the  following relations:
\[
\begin{split}
 Y_n(u) &= Z_n(u)  - 2n(2n-1 +\beta)(4n -3)Z_{n-1} + (2n)_3(2n-1 + \beta)_3Z_{n-2}\\
 uZ_n(u)  &= Y_{n+1} - (2n+1)(2n +\beta)(4n + 1)Y_{n} + (2n+1)_3(2n+ \beta)_3Y_{n-1}.
   \end{split}
\]  
 After simple manipulations we obtain relations of the form:
 
 \[
 \begin{split}
 uY_n(u) &= Y_{n+1}  +    \sum\limits_{j=0}^{3} \chi_j(n)Y_{n-j}\\
 uZ_n(u)  &= Z_{n+1} + \sum\limits_{j=0}^{3} \kappa_j(n)Z_{n-j}.  
 \end{split}
 \] 
 We see that the new polynomial systems are again VOP, satisfying 5-term recurrence. Moreover both of them possess the Bochner's property - the corresponding differential equations are easily derived from the differential equation for $P_n$. Let us write them explicitly  in the case $\beta = 1$ (for simplicity). First we transform the operator $L$ by the change of the variables $u = x^2$. We are going to use that     $\partial_x = 2 u^{1/2}\partial_u$ and $x\partial_x = 2u\partial_u$. Then the new operator is
 
  \begin{eqnarray*}
  	L_1 = [(2u\partial_u +1)2 u^{1/2}\partial_u]^2   + 2u\partial_u   
  	=  4[(2u\partial_u +1) u^{1/2}\partial_u]^2   + 2u\partial_u. 
  \end{eqnarray*}
 Simple computations show that

 \begin{eqnarray*}
 	L_1 	=  16(u\partial_u +1/2)^2 (u \partial_u  +1)\partial_u  + 2u\partial_u. 
 \end{eqnarray*}
 
If we put $L_Y = L_1/2$  the differential equation will read

\begin{eqnarray*}
	L_Y Y_n 	= nY_n.
\end{eqnarray*}
We see that    the polynomial system $\{Y_n\}$ corresponds to  $G = R(H)\partial_u$, where $R(H)= 8 (H+1)(H+1/2)^2$. 

The differential operator $L_3$, for which the VOP $\{Z_n\}$ are eigenfunctions, can be  easily obtained    having in mind that 
\[
L [u^{1/2}   Z_n(u^{1/2}) ]  =   n Z_n(u^{1/2}).
\]
   We see that $L_3$ can be obtained from the above $L_1$ by conjugation:

\[
L_3 = u^{-1/2}L_1 u^{1/2} =   16(u\partial_u +3/2)^2(u\partial_u)\partial_u + 2u\partial_u.
\]
Then we can put $L_Z = L_3/2$. We conclude that the VOP $\{Z_n\}$   come from $G$ with $G = 8(u\partial_u +3/2)^2(u\partial_u)\partial_u$. Notice that we need to change the variable $n$ into  $m = n +1/2$ to get the form of the present paper. 

Let us recall that 
 the classical Hermite polynomials   $H_n(x)$ have two different representation via Laguerre polynomials for the even and odd polynomials. Namely, 
 the even-indexed polynomials   $H_{2n}(x)$ are up to a multiplicative constant $L_n^{(-1/2)}(x^2)$, while the  odd-indexed  ones  $H_{2n+1}(x)$ are expressed up to a multiplicative constant as $xL_n^{(1/2)}(x^2)$.  We see the same phenomenon here: the even-indexed polynomial system $P_{2n}(x)$  are given by $Y_n(x^2)$, while    the odd-indexed ones  $P_{2n +1}(x)$ are expressed in terms of the polynomials $xZ_n(x^2)$.

 \eex

  \bre{ex3-11}
 \exref{ex2-2} can  be turned into a more general result. Consider any $G = R(H)\partial_x$,  and put $q(G) = G^l$. Then one can split the VOP into $l$ families $F_s = \{P_{lm+s}, \; m = 0, 1,2, \ldots \}, \; s = 0, 1, \ldots, l-1$. We leave the details to the reader. However to find the corresponding differential operators or which is the same - to find the  corresponding operators $G$ that generate the families requires some labor.   See \cite{Ho3} for more conceptual approach based on generalized hypergeometric functions.

 \ere
  
  \qed
 
 \bex{ex3-1}
 Let us take   the polynomial   $R(H) = H^2 +\al H + \be$, $G = R(H)  \partial_x$ and 
 $q(G) = G $. Then

 \[
 L_1 =  H + q'(G)G =  x^2\partial_x^3 + \alpha x\partial_x^2 + \beta \partial_x   + x\partial_x.  
 \]
 The corresponding polynomials satisfy the 4-terms recurrence relation 
 \begin{eqnarray*} 
 	xP^q_n(x)  &= &P^q_{n+1} -(3n^2+(2\alpha -1)n +\be)P^q_n(x) 
 	+ nR(n-1)(3n  + \al -2)P^q_{n-1}\\
 	&- &n(n-1)R(n-1) R(n-2) P^q_{n-2}.
 \end{eqnarray*}

 This example has been produced in \cite{BCD2} within the context of vector orthogonal polynomials which are generalized hypergeometric  functions.
 \eex

 \bex{ex3-2}
 Let us take   the polynomial 
 $q(G) = G^2/2 $,    $G =  (x\partial_x +1)^2\partial_x$. Then

 \[
 L_1 =  H + q'(G)G = [(x\partial_x +1)^2\partial_x]^2  + x\partial_x.
 \]
 This example seems to be new as  most of the examples originating from Sect. \ref{con} and  \ref{dis}.
 
 The  polynomials satisfy a  7-term recurrence relation which we skip as it is not much simpler than the corresponding general formula \eqref{new-rec3}. Here are the first few polynomials:
 \[
 P_0(x) = 1,  P_1(x) = x,  P_2(x) = x^2 +2,  P_3(x) = x^3 + 3.2x,  P_4(x) = x^4 + 4.3 x^2 +4!, \ldots
 \]
 It is clear that the arguments from \exref{ex2-2}  can be repeated. However the computations require more work. For this reason we again refer to \cite{Ho3}   for different approach.
 \eex

 \subsection{Discrete VOP}
 
Many  examples about continuous VOP can be easily reformulated into   examples about discrete VOP.  This can be done either by using the isomorphisms between the two classes or directly repeating the arguments.  
 
 \bex{ex1-d}
 
 We start with the algebra $\B_1$ generated by $\Delta, x\nabla, \hat{x}$. Take the simplest operator   $G(\Delta)  = - a\Delta$. We find that 
 
 \[
 \tilde{L} = a\Delta +   x \nabla.
 \]
 Define the polynomials  
 
 \[
 C_n   =e^{-a\Delta}(-x)_n = \sum_{j=0}^n \frac{(-a)^j\Delta^j(-x)_n}{j!}    =    \sum  \frac{(-a)^{j}(-n)_j(-x)_{n-j}}{j!}.
 \]
Simple computations (see \cite{Ho2}) show that
these are   the normed Charlier polynomials.  This example, together with the identical construction of Appell polynomials in \cite{Ho1} motivates  the name "Charlier-Appell polynomials" for  general  $q(\Delta)$.
 
 \eex

 \bex{M} As shown in  \cite{Ho2} the Meixner polynomials  can be obtained using the algebra   $\B_1$, generated by $\Delta,  G= (\Delta - \beta)x\nabla$  and $\hat{x}$.
 We defined $L$ starting from $(c-1)H = (1-c)x\nabla$ and $q(G) = \mu G$. Then  
 
 \[
 L = \mu(1-c) x(\Delta + \nabla) - \mu\beta\Delta + (1-c) x\nabla.
 \]
   If we  take   the constant  $\mu$  to be
 
 \[
 \mu = \frac{c}{1-c}.
 \]
 we obtain exactly  the Meixner operator
 
 \[
 L_1 = c (x+\beta)\Delta +   x\nabla.
 \]
 as given in \cite{KLS} for $\beta \neq -N, \; N \in \Nset$.     Thus the Meixner polynomials $M(x,\beta,c)$  coincide with our polynomial systems built from $G= (\Delta - \beta)x\nabla$.   See also \exref{L-M}.
 
 In the case when    $\beta = - N, \; N \in \Nset$ these are the Kravchuk polynomials $K_n(x, N, p)$ where $p = c/(c-1)$. 
 
 \eex
 
 \bex{M2}
 Here is the simplest new example. This time we start with $H = -x\nabla$  Let us take $q(G) =   G^2/2$, where  $G = (-x\nabla + \beta)\Delta$. Then the new polynomials  
 
 \[
 P_n(x) = \sum_{j=0}^{\infty} \frac{(G^2)^j (x)_n}{j!}
 \]
 are eigenfunctions of the operator $L = - x\nabla   + [(\nabla - \beta)x\Delta]^2$  
 
 \[
 L P_n(x) =   nP_n(x).
 \]

 The recurrence relation reads
 
 \begin{eqnarray*}
 	xP_n &=&  P_{n+1} +nP_n\\
   	 &=& -(n)_2(n-1 +\beta)^2P_{n-1} + (n)_4(n -1 + \beta)_4P_{n-4}.
 \end{eqnarray*}

Here are the first few polynomials:

\[
P_0=1, \; P_1 = x,\; P_2 = (x)_2 +2, \; P_3 = (x)_3 + 3x, \;  P_4 = (x)_4 +6(x)_2 + 3. 
\]

 One is tempted to repeat the construction of \exref{ex2-2} and build again two families of polynomials. Again one can take the even and   odd index polynomials. However it is not clear what can be the analog of the transformation $x^2 = u$.
 
 \eex

 \section{Connection between continuous and discrete VOP} \label{conn}

 We have seen   that the classical discrete and continuous systems are connected by a simple transformation. Let us examine some specific examples.
 \bex{L-M}
 
 From \exref{ex2-1} we know that    the Laguerre polynomials can be constructed using the operator $G^c = -(x\partial_x + \alpha +1)\partial_x$. 
 
 The Meixner polynomials $M_n(x;  \beta, c)$  are obtained  in \exref{M} from  $(x)_n$  via the operator  $G^d =  \mu G$, $G = (-x\nabla + \beta)\Delta$, where

 \[
 \mu =  \frac{c}{1-c}, 
 \]
 
 Let us consider the continuous   polynomials defined from   
 $x^n$ via $G^c = \mu (x\partial_x + \beta)\partial_x$ with an appropriate choice of the parameters $\mu, \beta$.  Simple computations show that up to the linear change of the variables $x =    -\mu t$ these are the Laguerre polynomials. 
 
 Then we can state the following  
 \bco{HL-CM}
 Meixner polynomials are the image of the Laguerre polynomials under the modified Mellin transform:
 \[
 \mathcal{M}^*(L^{\alpha}_n(-\mu t) ) = M_n(x; \beta, c),
 \]
 where 
 
 \[
 \mu = \frac{c}{1-c},\;\; \alpha = - \beta - 1.
 \]
 \proof Follows from     \prref{tran}. 
 
 \qed
 
 \eco

 \eex
 \bre{L-M}
 It is strange that such a simple formula has not been noticed till now.
 There are several papers where one can find maps between Laguerre and Meixner polynomials. In \cite{AA}  the authors also use Mellin transform and find similar formula but identify the final result as Meixner-Polaczek polynomials, which in fact are (from algebraic point of view) a specialization of Meixner polynomials.  See also \cite{Koo1} and \cite{Sz}.
 \ere
 Other examples of such dualities include Kravchuk polynomials and Laguerre polynomials with $\alpha = - N$. 
 
 \bex{H-mC}
 In this example we show the effect of the "base change" $x \rightarrow g = x  - H$. 
 Let us take $\B_1$,  generated by $x\partial_x, G= a\partial_x^2$ and $\hat{x}$. The corresponding polynomial system $e^G x^n$ consists of the Hermite polynomials up to rescaling $x \rightarrow \mu x$. The   discrete counterpart VOP becomes  
 
 \[
 H^d_n(x) = e^{a\Delta^2}(x)_n.
 \]
 
 However it is not orthogonal system but VOP with $d =2$. The explanation is that while the Hermite  polynomials satisfy the recurrence relation
 
 \[
 xH_n(x)  =  H_{n+1}  +nH_{n-1},
 \]
 their discrete counterparts   $H^d_n(x)$ satisfy

 \[
 gH^d_n(x) = H^d_{n+1}(x)  +n H^d_{n-1}(x).
 \]
 Using that   $g = x - H$   this  gives   $b'(x)      =   b'(g) + b'(H)$. We have 
 
 \[
 b'(H) = b(\sigma^{-1} (H)) =    b(H - 2a\Delta^2) = n -2a ( nT^{-1})^2.
 \]
  In explicit form the recurrence relation reads:

 \[
 xH^d_n(x) = H^d_{n+1}(x)  + n H^d_n(x)    +n H^d_{n-1}(x)  - 2an(n-1)H^d_{n-2} (x).
 \]

 \eex
 
 \bex{C}
 Let us  find the continuous counterpart  of Charlier polynomials. We start with $G = a\partial_x$. The corresponding automorphism is given by

 \[
 \sigma = e^{\ad_{a\partial_x}}.
 \]
 We are interested in the anti-involution $b'$ induced by $G$.  We have
 
 \[
 \begin{cases}
 b'(x) = T - a\\
 b'(H) = n + anT^{-1}\\
 b'(\partial_x) = T.
 \end{cases}
 \]
 In particular the recurrence relation is
 
 \[
 x\psi(x,n) = \psi(x,n+1) - a\psi(x,n)
 \]
 and the polynomials are simple shifts of the basic system:

 \[
 \psi(x,n) = e^{a\partial_x} x^n = (x + a)^n.
 \]
 At the same time their discrete counterpart - Charlier polynomials satisfy 
 
 \[
 gC_n(x) = C_{n+1}(x) - a C_n(x), 
 \]
 which, using $b'(g)   = b'(x- H)$ gives
 
 \[
 xC_n(x) = C_{n+1}(x) + (-a +  n) C_n(x) + an C_{n-1}(x).
 \]
 
 \eex

 

\end{document}